\documentclass[10pt, a4paper, margin=2.5cm]{article}

\usepackage[utf8]{inputenc}
\usepackage{amsmath} 
\usepackage{geometry}
\usepackage{amssymb} 
\usepackage{amsthm}
\usepackage[T1]{fontenc}
\usepackage{tikz-cd}
\usepackage[english]{babel}

\newtheorem*{qu}{Question}

\newtheorem{proposition}{Proposition}[section]
\newtheorem{Remark}[proposition]{Remark}

\newtheorem{fact}[proposition]{Fact}
\newtheorem{definition}[proposition]{Definition}
\newtheorem{theorem}[proposition]{Theorem}
\newtheorem{Lemma}[proposition]{Lemma}
\newtheorem{corollary}[proposition]{Corollary}

\newcommand{\ad}{\operatorname{ad}}

\newcommand{\im}{\operatorname{im}}
\newcommand{\End}{\operatorname{End}}

\newcommand{\ann}{\operatorname{ann}}
\newcommand{\Aut}{\operatorname{Aut}}
\newcommand{\Hom}{\operatorname{Hom}}
\newcommand{\Der}{\operatorname{Der}}
\newcommand{\IDer}{\operatorname{IDer}}
\newcommand{\res}{\operatorname{res}}

\title{On the cohomology of finite-dimensional nilpotent groups and Lie rings\footnote{ Keywords: nilpotent group, nilpotent Lie algebra, cohomology, finite-dimensional theory
\\
MSC 2020 : 20J06 ,  03C60   }}

\author{Samuel Zamour}

\begin{document}
\maketitle
\begin{abstract}
We establish vanishing results for the first cohomology group of nilpotent groups and Lie rings when the submodule of invariants is trivial. Our results are obtained within a model-theoretic setting, namely for structures that are definable in a finite-dimensional theory, which encompasses algebraic groups over algebraically closed fields, real semi-algebraic groups, and finite-dimensional Lie algebras over  an algebraically or real closed field. Since classical tools - such as computations with spectral sequences and rigidity of the linear dimension - are not available in our setting, we develop an elementary algebraic approach. As applications, we derive a form of Frattini's argument for Cartan subrings and a definable version of Maschke's theorem for actions of definable connected $p$-divisible abelian groups, with a view toward the ongoing study of soluble finite-dimensional Lie rings.
\end{abstract}

\section{Introduction}
The present paper lies at the interface of algebra and model theory. 
Its primary aim is to develop cohomological tools for the study of groups and Lie rings in the model-theoretic setting of \emph{finite-dimensional theories} \cite{Wagd}, a framework we now describe informally; precise definitions are given in Section 2.

A \emph{finite-dimensional} structure is a structure equipped with a well-behaved notion of dimension on its definable sets, analogous to the Zariski dimension in algebraic geometry but requiring no underlying field or topology. From an algebraic point of view, this framework encompasses algebraic groups over algebraically closed fields, semi-algebraic groups over the reals, and Lie algebras of finite linear dimension over algebraically closed or real closed fields. One key feature of this setting is that one can work with a meaningful notion of \emph{length} of a module (defined via chains of definable connected submodules) and one can use Schur-like linearization principles, while having no access to any ambient geometry or topology. In this context, a \emph{Lie ring} is an abelian group equipped with a bracket satisfying the axioms of a Lie bracket--- antisymmetry, bi-additivity, Jacobi identity --- but living in no ambient vector space; it is the correct analogue of a Lie algebra when no linear structure is present from the outset.

\subsection*{Background and motivation}

Cohomological methods have long been central to the structural theory 
of both groups and Lie algebras. In the Lie-theoretic setting, Barnes 
\cite{Barn1, Barn2} used them to obtain fundamental results for 
soluble Lie algebras of \emph{finite linear dimension}. His starting 
point \cite{Barn1} is the following vanishing theorem:

\begin{fact}[Barnes]\label{fact:barnes1}
Let $\mathfrak{g}$ be a nilpotent Lie algebra of finite linear 
dimension over a field, and let $A$ be a $\mathfrak{g}$-module. 
Suppose $A^{\mathfrak{g}} = H^0(\mathfrak{g},A) = 0$. Then:
\begin{itemize}
  \item $H^n(\mathfrak{g},A)=0$ for all $n \geq 0$;
  \item $(B/C)^{\mathfrak{g}}=0$ for all submodules $C \subseteq B 
  \subseteq A$.
\end{itemize}
\end{fact}

Analogous vanishing results hold for nilpotent finite groups and nilpotent algebraic groups, though the proofs differ in each case. Building on Fact~\ref{fact:barnes1}, Barnes derived 
the following structural consequence, a result in the spirit of Frattini's argument \cite{Barn2}:

\begin{fact}[Barnes]\label{fact:barnes2}
Let $\mathfrak{g}$ be a Lie algebra of finite linear dimension, and 
let $\mathfrak{c}$ be a Cartan subalgebra of an ideal 
$\mathfrak{i} \trianglelefteq \mathfrak{g}$. Then
\[
\mathfrak{g} = \mathfrak{i} + N_{\mathfrak{g}}(\mathfrak{c}).
\]
\end{fact}

These results are classical when $\mathfrak{g}$ is a Lie algebra of finite linear dimension. However, if one tries to carry them over to the model-theoretic setting of finite-dimensional theories the classical proofs break down completely. The obstruction is genuine: the standard proof of Fact~\ref{fact:barnes1} proceeds by 
induction on $\dim_k \mathfrak{g}$ and by exploiting the 
Hochschild--Serre spectral sequence $H^p(\mathfrak{g}/\mathfrak{h}, H^q(\mathfrak{h},A)) \Rightarrow H^{p+q}(\mathfrak{g},A)$ for an ideal $\mathfrak{h}$, neither of which are available when there is no underlying field $k$ or linear dimension. In the context of finite-dimensional theories, the notion of length of a definable connected module is well-defined; however, it does not enjoy the rigidity of linear dimension (or cardinality in the case of finite groups) and behaves quite differently in inductive arguments. 

We must therefore return to elementary cochain computations and exploit the abstract dimension-theoretic properties of definable sets to replace geometric or topological arguments.

\subsection*{Main results}

We prove the following generalisation of Facts~\ref{fact:barnes1} and~\ref{fact:barnes2} to the model-theoretic setting of finite-dimensional theories. We first formulate our main cohomological result for groups : group cohomology provides a standard and flexible framework and our approach yields a uniform treatment encompassing both nilpotent finite groups and nilpotent algebraic groups, in the spirit of Barnes’ methods for Lie algebras. The corresponding statement for Lie rings is established in Section~3.

\begin{theorem}\label{thm:main1}
In a finite-dimensional$^\circ$ theory, let $G$ be a definable connected nilpotent group, and let $A$ be a definable connected $G$-module. Suppose $A^G = 0$. Then:
\begin{enumerate}
\item $H^1(G,A) = 0$.
\item For all definable connected submodules $V \leq U \leq A$, one has $(U/V)^G = 0$.
\end{enumerate}
\end{theorem}

As algebraic applications, we establish two structural results.

\begin{theorem}\label{thm:main2}
\begin{enumerate}
\item \textup{(Definable Maschke's theorem)} In a 
finite-dimensional$^\circ$ theory, let $T$ be a definable connected abelian group with trivial $p$-torsion, and let $V$ be a definable connected $T$-module $p$-elementary as an abelian group. Suppose $V^T = 0$. Then
\[
V = V_1 \oplus \cdots \oplus V_r,
\]
where each $V_i$ is a definable connected irreducible $T$-module.

\item \textup{(Frattini's argument for Lie rings)} In a finite-dimensional$^\circ$ theory, let $\mathfrak{g}$ be a definable connected Lie ring, and let $\mathfrak{c}$ be a Cartan subring of a definable connected ideal 
$\mathfrak{i} \trianglelefteq \mathfrak{g}$. Then
\[
\mathfrak{g} = \mathfrak{i} + N^{\circ}_{\mathfrak{g}}(\mathfrak{c}).
\]
\end{enumerate}
\end{theorem}

Although Theorem~\ref{thm:main1} concerns only the first cohomology group, it already has significant structural consequences. In particular, as demonstrated in \cite{TZ}, the Frattini-type result of Theorem~\ref{thm:main2}(2) underlies a nascent Frattini theory for 
soluble finite-dimensional Lie rings, in which Cartan subrings play the role of Cartan subalgebras in classical soluble Lie algebra theory. Even when restricted to the case of groups of finite Morley rank --- one of the most studied instances of finite-dimensional theories --- 
our cohomological approach to these results appears to be new.

\subsection*{Context and applications}

The structural theory of soluble Lie rings is a key ingredient in the program of classifying simple infinite Lie rings of finite Morley rank. This program is the Lie-theoretic analogue of the celebrated Cherlin--Zilber conjecture on groups of finite Morley rank. It proposes that a simple infinite Lie ring of sufficiently large positive characteristic and finite Morley rank is isomorphic to a simple Lie algebra of finite linear dimension over an algebraically closed field (see \cite{DT1,DT2} for the current state of the program). The analogy with the Cherlin--Zilber conjecture motivates the development of a purely algebraic structural theory of soluble finite-dimensional Lie rings --- one in 
which the geometric methods of classical Lie theory (Killing form, root systems, representation theory over $\mathbb{C}$) are simply unavailable and must be replaced by elementary algebraic arguments. 
The results of the present paper are intended as basic tools in this direction.

\subsection*{Methods and limitations}

A word is in order on the reason why we restrict to $H^1$. The definable cohomology, $H^n_{\mathrm{def}}(G,A)$, in which cochains are required to be definable, behaves poorly in the absence of additional geometric structure: in particular, short exact sequences of $G$-modules need 
not give rise to long exact sequences in definable cohomology. We therefore work with the classical (non-definable) cohomology, which does admit such long exact sequences. Our two main algebraic tools are: 
(i) the exactness of the inflation-restriction sequence in degree~$1$, 
proved by elementary cochain computations; (ii) the long exact cohomology sequence associated to a short exact sequence of $G$-modules. These suffice for the proofs of Theorems~\ref{thm:main1} 
and~\ref{thm:main2}. The Hochschild--Serre spectral sequence plays no role.

By contrast, establishing the vanishing of $H^n(G,A)$ for $n \geq 2$ under the hypotheses of Theorem~\ref{thm:main1} appears to be significantly harder. The difficulty is not merely technical: in the model-theoretic setting of finite-dimensional theories, the standard dimension arguments which drive inductive proofs of higher vanishing (such as 
those of Barnes \cite{Barn1}) are not available, and it is unclear whether the conclusion can hold in this generality. We leave the following as an open question:

\begin{qu}\label{qu:higher}
In a finite-dimensional theory, let $G$ be a definable connected nilpotent group, and let $A$ be a definable connected $G$-module with $A^G = 0$. Does $H^n(G,A)$ vanish for all $n \geq 1$? Does the analogous statement hold for the definable cohomology groups $H^n_{\mathrm{def}}(G,A)$?
\end{qu}

A positive answer to this question would provide, in particular, a purely model-theoretic generalisation of Fact~\ref{fact:barnes1} in full. We remark that even a negative answer --- a counterexample to higher vanishing in a finite-dimensional theory 
--- would be of considerable interest, as it would delineate a genuine difference between the model-theoretic and the linear-algebraic 
settings.

\subsection*{Organisation of the paper}

Section 2 recalls the necessary background on group cohomology, introduces the notion of finite-dimensional$^\circ$ theory, and proves Theorem~\ref{thm:main1} together with the definable Maschke theorem (Theorem~\ref{thm:main2}(1)). Section 3 develops the parallel theory for Lie rings and proves 
Theorem~\ref{thm:main2}(2).

\medskip
\noindent\textbf{Acknowledgements.} We thank Adrien Deloro for his many comments and suggestions. We are also grateful to Moreno Invitti for his key observation regarding long exact cohomological sequences. We are also grateful to the anonymous referee for many valuable comments that have improved a first version of this paper.
\section{Groups}
\subsection{Preliminaries on the cohomology of groups}
We recall some basic facts regarding the cohomology of groups. Although the results presented in this subsection are fairly standard, most textbooks adopt a categorical approach that does not extend well to abstract Lie rings, since general cohomological results may fail for higher cohomology groups in that setting. Moreover, it obscures the concrete interpretation through cochains which plays a central role in our approach. We therefore provide elementary algebraic proofs that extend naturally to the case of Lie rings. For general references on the topic of the cohomology of groups, we refer to \cite[Chap. 6]{Weib} and to \cite[Chap. VII]{Ser} (the first one is categorical in nature while the second one insists more on cochains).
\begin{definition}
A $G$-module is the data of ($G$,$A$), where $G$ is a  group and $A$ is an abelian group, equipped with a group action "$\cdot$" of $G$ on $A$. Equivalently, this amounts to giving a group homomorphism $\rho$ from $G$ to $\Aut(A)\subseteq \End(A)$.
\end{definition}

\begin{definition}
  Let $(G,A)$ be a $G$-module. We define a cochain complex as follows:
\begin{itemize}
    \item $C^0(G,A)=A$.
    \item For every integer $n\geq 1$, $C^{n}(G,A)$ denotes the abelian group of functions from $G^n$ to $A$.

    \end{itemize}
    The differential $d_n$ from $C^n(G,A)$ to $C^{n+1}(G,A)$ is defined, for every map $f$ in $C^n(G,A)$, by:

\[\begin{aligned}& d_nf(g_1,\dots,g_{n+1})\\
& =g_1\cdot f(g_2,\dots,g_{n+1})+(-1)^{n+1}f(g_1,\dots,g_n)\\
& + \sum_{1\leq s \leq n}(-1)^s f(g_1,\dots,g_sg_{s+1},\dots,g_{n+1})\end{aligned}\]
\end{definition}
We have $d_{n+1}\circ d_n=0$. We set $C(G,A)=(C^n(G,A))_{n\in\mathbb{N}}$.
\\
\\
 We may now define the associated cohomology groups:
 \begin{definition}
  Given the complex defined above, we introduce for every positive integer $n$ the $n$-th  cohomology group, denoted by $H^n(G,A)$ : 
  \[H^n(G,A)=\ker(d_n)/\im(d_{n-1}).\]
 \end{definition}
 The group $G$ acts on the cochain groups of the form $C^n(H,A)$, for $H$ a  normal subgroup, in a way that is compatible with the differentials : 
 \\
 let $f\in C^n(H,A)$ and $x\in G$, then \[(x\cdot f)(h_1,\dots,h_n)=x\cdot f(h_1^x,\dots,h_n^x).\]
 \\
 \\
 For $n=0$ and $n=1$, it is possible to give a convenient algebraic description of the cohomology groups.
\begin{fact}\label{carac deriv}
  \begin{enumerate}
     \item $H^0(G,A)=A^G=\{a\in A : g\cdot a=a, \text{for all $g\in G$}\}$.
     \item $H^1(G,A)=\Der(G,A)/\IDer(G,A)$, where \[\Der(G,A)=\{f\in C^1(G,A)  : f(xy)=x\cdot f(y)+f(x)\}\] and \[\IDer(G,A)=\{f\in C^1(G,A) : f(x)=x\cdot a-a~ \text{for some $a$ in $A$}\}.\]
     In particular, if $f\in \Der(G,A)$, then $f(1)=0$.
 \end{enumerate}  
\end{fact}
\begin{Remark}\label{trivial action}
We observe that $G$ acts trivially on $H^1(G,A)$: let $[f]\in H^1(G,A)$, where $f$ is a derivation, and let $x\in G$; then \[\begin{aligned}&(x\cdot f)(y)-f(y)=x\cdot f(y^x)-f(y)\\
&=(f(yx)-f(x))-f(y)=y\cdot f(x)-f(x)\end{aligned}\] for all $y\in G$; in other words, $(x\cdot f)-f$ is an inner derivation and thus $[x\cdot f]=[f]$.
 \end{Remark}
 We now turn to the question of the relations between exact sequences of $G$-modules and long exact sequence of cohomology groups.

\begin{fact}
    Every short exact sequence of $G$-modules $0\rightarrow A\overset{i}{\rightarrow}B\overset{p}{\rightarrow}C\rightarrow0$ induces the following short exact sequence : 
    \[0\rightarrow C^1(G,A)\overset{i^*}{\rightarrow} C^1(G,B)\overset{p^*}{\rightarrow} C^1(G,C),\] 
    where $i^*,p^*$ are defined by $i^*(f)=i\circ f$ and $p^*(f)=p\circ f$. Moreover, $i^*$ and $p^*$ commute with the differential. 
\end{fact}
\begin{proof}
    Since $i$ and $p$ commute with the action of $G$, the formula for the differential allows us to show that $i^*,p^*$ commute with it.
    
     Let us show that $i^*$ is injective. Indeed, for $\phi \in C^1(G, A)$, $i^*(\phi)=0$ iff $i(\phi(x))=0$ for all $x\in G$ iff $\phi(x)=0$ for all $x\in G$ by injectivity of $i$.

    Now, we turn to the equality $\im(i^*)=\ker(p^*)$. Let $\phi\in \im(i^*)\subseteq C^1(G, B)$; then there exists $\psi\in  C^1(G, A)$, such that $\phi=i\circ \psi$ and hence $p^*(\phi)=p\circ i \circ \psi$. Consequently, for every $x\in G$, $p(i(\psi(x))=0$ since $\im(i)\subseteq \ker(p)$, and thus $\phi\in \ker(p^*)$.
    
    Conversely, let $\phi\in \ker(p^*)$; then $p(\phi(x))=0$ for every $x\in G$. Thus $\phi(x)\in \im(i)$, since $\im(i)=\ker(p)$. Therefore for every $x$ there exists an unique $y_{x}\in A$, such that $i(y_{x})=\phi(x)$ by injectivity of $i$. Now, it suffices to consider the following element of $C^1(G,A)$ : $\psi :x\mapsto y_{x}$.
\end{proof}
\begin{corollary}
    Every short exact sequence of $G$-modules $0\rightarrow A\overset{i}{\rightarrow}B\overset{p}{\rightarrow}C\rightarrow0$ induces the following short exact sequence :
    \[0\rightarrow \Der(G,A)\overset{i^*}{\rightarrow} \Der(G,B)\overset{p^*}{\rightarrow} \Der(G,C),\]
    where $i^*,p^*$ are defined by $i^*(f)=i\circ f$ and $p^*(f)=p\circ f$.  
\end{corollary}
\begin{proof}
    Since $i$ and $p$ are $G$-morphisms, the sequence is well-defined. Moreover, let $\phi\in \ker(p^*)$ and consider the map $\psi: G\rightarrow A,  x\mapsto y_{x}$, where $i(y_{x})=\phi(x)$. Then, for all $x,y\in G$, we have $i(x\cdot\psi(y)+\psi(x))=x\cdot i(\psi(y))+i(\psi(x))=x\cdot\phi(y)+\phi(x)=\phi(xy)=i(\psi(xy))$, so, by injectivity of $i$, $\psi(xy)=x\cdot\psi(y)+\psi(x)$.
\end{proof}
One important feature of cohomology is the possibility of getting a "long" exact sequence of cohomology groups from a short exact sequence of $G$-modules.
\begin{fact}\label{long cohomology sequence}

    Let $0\rightarrow A_0\rightarrow A_1 \rightarrow A_3\rightarrow 0$ be a short exact sequence of $G$-modules. Then there exists a "long" exact sequence of cohomology groups of the following form :
    \[\begin{aligned}&A_0^G\rightarrow A_1^G\rightarrow A_2^G\\
    &\overset{\delta_0}{\rightarrow}H^{1}(G,A_0)\rightarrow H^{1}(G,A_1)\rightarrow\ H^{1}(G,A_2)\end{aligned}\]
The morphism $(\delta_0)$ is called a connecting morphism.
\end{fact}
\begin{proof}
    By Fact \ref{carac deriv}, it suffices to apply the Snake Lemma \cite[Lemma 9.1, Chap. III]{Lan} to the following commutative diagram, where $\phi_i:A_i\rightarrow \Der(G,A_i), a_i\mapsto(g\mapsto g\cdot a_i-a_i), i=0,1,2 $ :
\[
\begin{tikzcd}
{} &
A_0 \arrow[r,"i"] \arrow[d,"\phi_1"] &
A_1 \arrow[r,"p"] \arrow[d,"\phi_2"] &
A_2 \arrow[r] \arrow[d,"\phi_3"] &
0
\\
0 \arrow[r] &
\mathrm{Der}(G,A_0) \arrow[r,"i^{*}"] &
\mathrm{Der}(G,A_1) \arrow[r,"p^{*}"] &
\mathrm{Der}(G,A_2)
\end{tikzcd}
\]
\end{proof}
We now introduce the Hochschild-Serre inflation-restriction sequence. 
\begin{definition} Let $(G,A)$ be a $G$-module and let $H$ be a normal subgroup of $G$.
    \begin{itemize}
        \item For every $n\in \mathbb{N}$, we define the restriction map $\res : C^n(G,A)\longrightarrow C^n(H,A), f\mapsto f|_{H}$.
        \item We define the inflation map $\inf : C^n(G/H, A^H)\longrightarrow C^n(G,A), f\mapsto f\circ \pi$, where $\pi$ denotes the canonical projection from $G$ to $G/H$. 
    \end{itemize}
\end{definition}
It would have made sense to define the inflation map on $C^n(G/H,A)$ but exactness would have been lost in the subsequent Fact \ref{suite inflation-restriction groupe}.
Since the inflation and restriction morphisms commute with the differentials, we can pass to the cohomology groups and we get a sequence $H^1(G/H,A^{H})\overset{\inf}{\rightarrow} H^1(G,A)\overset{\res}{\rightarrow} H^1(H,A)$. In fact, this sequence is exact, and it can be proved by an elementary algebraic argument (see \cite[Proposition 4, §6, Chap. VII]{Ser} for a proof using computations on cochains).
\begin{fact}\label{suite inflation-restriction groupe}
The sequence \[0\rightarrow H^1(G/H,A^{H})\overset{\inf}{\rightarrow} H^1(G,A)\overset{\res}{\rightarrow} H^1(H,A)\] is exact.  
\end{fact}
Note that we slightly abuse notation by using $\inf$ and $\res$ instead of $\inf^*$ and $\res^*$.

We can refine the image of the restriction morphism.
\begin{fact}\label{image restriction inflation groupe}
    $H^1(G,A)\overset{\res}{\rightarrow}H^1(H,A)^{G/H}$.
\end{fact}
\begin{proof}
    Let $[\psi] \in H^1(H,A)$ such that $[\psi]=[\phi|_{H}]$ for $\phi\in \Der(G,A)$. Since $\phi$ is a $G$-derivation, for all $x\in G$ and $h\in H$, we have $(x\cdot \psi)(h)-\psi(h)=x\cdot \phi(h^x)-\phi(h)=(\phi(hx)-\phi(x))-\phi(h)=h\cdot \phi(x)-\phi(x)$, and so $[x\cdot \psi]=[\psi]$ in $H^1(H,A)$.
\end{proof}
The following lemma will be used repeatedly.
\begin{Lemma}\label{fidélité groupe}
    Let $(G,A)$ be a  $G$-module and let $H$ be a  central subgroup of $G$, such that $H\leq C_{G}(A)$. Suppose that $A^{G}=0$. Then $H^1(H,A)^{G/H}=0$ and $H^1(G/H, A)\simeq H^1(G, A)$.
\end{Lemma}
\begin{proof}
Since $A^{H}=A$, by Fact \ref{suite inflation-restriction groupe} and Fact \ref{image restriction inflation groupe}, the sequence $0\rightarrow H^1(G/H,A)\overset{\inf}{\rightarrow} H^1(G,A)\overset{\res}{\rightarrow} H^1(H,A)^{G/H}$ is exact. Therefore, it suffices to prove that \[H^1(H,A)^{G/H}=0.\]
Let $f\in \Der(H,A)$; suppose that for all $x\in G$, $[x\cdot f]=[f]$. For all $x\in G$ there exists an $a_x\in A$ such for all $h\in H$ one has $x\cdot f(h^x)-f(h)=x\cdot f(h)-f(h)=h\cdot a_x-a_x=0$ since $H$ acts trivially on $A$ and $H\leq Z(G)$. In other words, for all $x\in G$ and $h\in H$, the equality $x\cdot f(h)=f(h)$ holds, and thus for all $h\in H$, $f(h)\in A^G=0$; finally, we obtain that $f=0$. 
\end{proof}
\subsection{Finite-dimensional nilpotent groups}
We will study group cohomology inside the framework of finite-dimensional theories.
\begin{definition} \cite[Definition 1.1]{Wagd}
    A theory $T$ is said to be finite-dimensional if there exists a dimension function, $\dim$, defined on the collection of all interpretable sets in models of $T$, with values in $\mathbb{N}$, such that for a formula $\phi(x,y)$ and interpretable sets $X, Y$, the following properties are satisfied :
    \begin{itemize}
        \item (Invariance) If $a\equiv a'$, i.e., $a$ and $a'$ have same types, then $\dim(\phi(x,a))=\dim(\phi(x,a'))$.
        \item (Finiteness) $\dim(X)=0$ iff $X$ is finite.
        \item (Union) $\dim(X\cup Y)=\max\{\dim(X),\dim(Y)\}$.
        \item (Fibration) If $f:X\longrightarrow Y$ is an interpretable map, such that $\dim(f^{-1}(y))\geq d$ for all $y\in Y$, then $\dim(X)\geq \dim(Y)+d$.
        \item (Low fibration) If $f:X\longrightarrow Y$ is an interpretable map, such that $\dim(f^{-1}(y))\leq d$ for all $y\in Y$, then $\dim(X)\leq \dim(Y)+d$.
    \end{itemize}
\end{definition}
Here, \textit{interpretable} has a technical meaning that ensures that a dimension can be assigned to quotients : a subset $X$ is said to be interpretable in a model $M$ of $T$ if there exist a definable subset $X'\subseteq M^n$ for some $n\in \mathbb{N}$ and a definable surjective map $f$ from $X'$ to $X$ such that $\{(\overline{x},\overline{y})\in M^{2n} : f(\overline{x})=f(\overline{y})\}$ is definable.

The dimension extends to type-definable sets, \textit{i.e.}, bounded intersections of definable sets, by taking the infimum of the dimensions of the definable sets which appear in the intersection. Here, \textit{bounded} means \textit{not increasing in elementary extensions} (and so \textit{unbounded} means \textit{increasing in elementary extensions}). In the following, we will abuse notation and we say definable instead of interpretable.

We will need to slightly strengthen our notion of finite-dimensionality. 
\begin{definition}
 A finite-dimensional theory is said to be finite-dimensional$^{\circ}$ if it satisfies the following property: 
 \\
 Any definable group admits a connected component, i.e., a smallest definable subgroup of finite index.
\end{definition}
The theory of a group of finite Morley rank and of a group definable in an o-minimal structure are examples of finite-dimensional$^{\circ}$ theories. This model-theoretic framework is suitable for inductive reasoning when working with groups.
\begin{fact}\cite[Proposition 2.3]{Wagd}
  Let $H\leq G$ be two type-definable groups in a finite-dimensional theory. If $H$ has infinite index in $G$, then $\dim(H)<\dim(G)$.  
\end{fact}
A definable group is said to be \textit{connected} if it has no proper definable subgroup of finite index.
\\
Another remarkable property of finite-dimensional theories concerns the possibility of obtaining a definable version of Schur's Lemma.
\begin{definition}
 In a finite-dimensional theory, let $V$ be a definable connected abelian group. We define $\End_{\text{def}}(V)$ as the ring of definable endomorphisms of $V$. Let $S$ be a subring of $\End_{\text{def}}(V)$; then $V$ is said to be an irreducible $S$-module if every proper definable connected $S$-invariant subgroup is trivial.
\end{definition}
\begin{fact}\label{linearisation}(adapted from \cite{DelZ})
    Let $T$ be a finite-dimensional theory. Let $V$ be a definable connected abelian group and $S\leq \End_{\text{def}}(V)$ an invariant, \textit{i.e.}, a bounded union of type-definable sets, and unbounded subring of definable endomorphisms. Suppose that $V$ is an irreducible $S$-module. Then $T=C_{\End_{\text{def}}(V)}(S)$ is a definable skew-field.
    \end{fact}
In the following, unless explicitly stated otherwise, an irreducible module will be understood in the sense indicated above, i.e., in reference to the definable connected category. 
\begin{definition}
In a finite-dimensional$^{\circ}$ theory, let $(G,A)$ be a definable $G$-module. Suppose that $A$ is connected. We define the $G$-length of $A$, $\lg_G$, as the maximum length of a composition series of the following form:
\[A_0=A\geq\ A_1\geq \dots \geq A_n=0,\] where each $A_i$ is a definable connected $G$-module such that $A_{i}/A_{i+1}$ is $G$-irreducible (in the definable connected category) for $i\in\{0,\dots,n-1\}$.
\end{definition}

\begin{Lemma}\label{abelian case finite-dimensional}
  In a finite-dimensional$^{\circ}$ theory, let $(G,A)$ be a definable $G$-module. Suppose that $G$ is connected and abelian, that $A$ is connected and that $A^G=0$. Then $H^1(G,A)=0$.   
\end{Lemma}
\begin{proof}
 First, according to Lemma \ref{fidélité groupe}, up to considering $G/C_{G}(A)$, we may assume that $G$ acts faithfully. The result is trivial for $A=0$. 
 We proceed by induction on the length, $\lg_G$, of $A$ as a $G$-module (in the definable connected category).
 \begin{itemize}
     \item First, suppose $\lg_G(A)=1$, \textit{i.e.}, $A$ is irreducible. We denote by $\rho : G\rightarrow \Aut_{\text{def}}(A)\subseteq \End_{\text{def}}(A)$ the action. By Fact \ref{linearisation}, we can linearize: there exists a definable skew-field $\mathbb{K}$ such that $\rho(G)\leq \mathbb{K}$. Up to considering $Z(\mathbb{K})$, which is also definable infinite, we may assume that $\mathbb{K}$ is commutative definable infinite. Now, let $f : G\rightarrow A$ be a non-trivial $G$-derivation. Let $y\in G-\{1\}$, such that $f(y)\neq 0$; then $(\rho(y)-Id)$ is invertible. Let $a=(\rho(y)-Id)^{-1}\cdot f(y)$. For $x\in G$, by commutativity of $G$ and the definition of a derivation, we have $(\rho(y)-Id)\cdot f(x)=\rho(x)\cdot f(y)-f(y)$, thus, \[f(x)=(\rho(y)-Id)^{-1}\cdot \rho(x)\cdot f(y)-(\rho(y)-Id)^{-1}\cdot f(y)=\rho(x)\cdot a-a\]
        using commutativity once more. This shows that $H^1(G,A)=0$.
     \item We may therefore assume $\lg_G(A)>1$; then there exists a non-trivial definable connected $G$-module $B$ strictly contained in $A$ such that $\lg_G(B)<\lg_G(A)$ and $\lg_G(A/B)<\lg_G(A)$.  Note that by hypothesis $B^{G}=0$. By induction on the length, we may therefore assume that $H^1(G,B)=0$. The long exact cohomology sequence will now come into play. Indeed, we have a short exact sequence of $G$-modules $0\rightarrow B\rightarrow A \rightarrow A/B \rightarrow 0$ which, by virtue of Fact \ref{long cohomology sequence}, induces the following long exact cohomology sequence : \[\begin{aligned}&0\rightarrow B^G\rightarrow A^G\rightarrow (A/B)^G\\
    &\overset{\delta_0}{\rightarrow}H^{1}(G,B)\rightarrow H^{1}(G,A)\rightarrow H^1(G,A/B)\end{aligned}\]
By hypothesis, $A^G=B^G=0$ and by induction hypothesis, $H^1(G,B)=0$, so the exactness of the sequence then yields $(A/B)^G=0$; we may therefore apply once more the induction hypothesis and use the exactness of the long cohomology sequence to have $H^1(G,A)=0$.
 \end{itemize} 
\end{proof}
Now, we pass to the nilpotent case. Recall that an infinite definable nilpotent group has infinite center (see \cite[Lemma 5.1, I]{ABC}).
\begin{theorem}\label{g-module}
In a finite-dimensional$^{\circ}$ theory, let $(G,A)$ be a definable $G$-module. Suppose that $G$ is connected and nilpotent, that $A$ is connected, and that $A^G=0$. Then $H^1(G,A)=0$.   
\end{theorem}
\begin{proof}
 We proceed by induction on the nilpotency class of $G$. By Lemma \ref{abelian case finite-dimensional}, we may suppose that $G$ is not abelian. Let $Z=Z(G)^{\circ}$; note that $G/Z$ and $Z$ are infinite definable connected nilpotent groups of nilpotency class strictly inferior.
 \\
First, suppose that $A$ is irreducible (in the definable connected category) and consider the definable submodule $A^{Z}$. Then we have two cases :
\begin{itemize}
    \item $A^{Z}=A$. By Lemma \ref{fidélité groupe}, $H^1(G,A)\simeq H^1(G/Z,A)$, and we are done by induction on the nilpotency class.
    \item $A^{Z}$ is finite. In fact, $A^{Z}=0$ : for $a\in A^Z$, one has $G\cdot a \le A^Z$, so $G\cdot a$ is both finite and connected, hence trivial. Notice that $H^1(G/Z,A^{Z})=0=H^1(Z,A)$ by Lemma \ref{abelian case finite-dimensional} and we are done by Fact \ref{suite inflation-restriction groupe}.
\end{itemize}
Now, we proceed by induction on the length of $A$ as a $G$-module exactly in the same way as in the proof of Lemma \ref{abelian case finite-dimensional}.
   \end{proof}

We conclude this subsection with a corollary on composition factors.
\begin{corollary}\label{facteur de composition}
In a finite-dimensional theory$^{\circ}$, let $G$ be a definable connected nilpotent group and let $A$ be a definable connected $G$-module. Suppose that $A^G=0$. Then for all definable connected submodules $V\leq U$ contained in $A$, we have $(U/V)^G=0$.
\end{corollary}
\begin{proof}
   Indeed, we observe that $U^G=0=V^G$; then it suffices to consider the long exact cohomology sequence : \[\begin{aligned}&0\rightarrow V^G\rightarrow U^G\rightarrow (U/V)^G\\
    &\overset{\delta_0}{\rightarrow}H^{1}(G,V)\end{aligned}\]  
    Theorem \ref{g-module} allows us to conclude.
\end{proof}

\subsection{Definable Maschke's theorem}

In this subsection, we provide an application of the cohomological study of finite-dimensional nilpotent groups by proving a definable version of Maschke's theorem. We closely follow the argument for the case of groups of finite Morley rank given by Tindzhogo Ntsiri in \cite{Tind}. Nevertheless, even for groups of finite Morley rank, the result presented here is new since we can weaken the hypothesis "decent torus" and replace it by "definable connected abelian p-divisible". We will need the following definability lemma:
\begin{fact}\label{field definability}(\cite[Proposition 3.6]{Wagd}, extracted from \cite{DelZ})
In a finite-dimensional theory, let $\mathbb{K}$ be a skew-field such that :
\begin{itemize}
    \item there is an upper bound on dimensions of type-definable subsets of $\mathbb{K}$.
    \item $\mathbb{K}$ contains an invariant unbounded subset.
\end{itemize}
Then $\mathbb{K}$ is definable.
\end{fact}
\begin{proposition}\label{anneau invariant}
In a finite-dimensional$^{\circ}$ theory, let $V$ be a definable connected $p$-elementary abelian group and let $R$ be an unbounded invariant commutative ring of definable endomorphisms of $V$. Suppose that $R^p=R$. Let $W$ be a definable connected irreducible $R$-submodule, such that $V/W$ is also $R$-irreducible. Suppose also that $R/\ann_R(W)$ and $R/\ann_R(V/W)$ are unbounded. Then $W$ admits a definable connected $R$-invariant complement.
\end{proposition}
\begin{proof}
By considering the action of $R$ on $W$ and $V/W$ and by applying Fact \ref{linearisation}, the rings $C_{\End_{\text{def}}(W)}(R/\ann_R(W))=K_W$ and $C_{\End_{\text{def}}(V/W)}(R/\ann_R(V/W))=K_{V/W}$ are in fact definable skew-fields. Note that $R/\ann_R(W)\hookrightarrow K_W$ and $R/\ann_R(V/W)\hookrightarrow K_{V/W}$; hence $\ann_R(W)$ and $\ann_R(V/W)$ are prime ideals. Moreover, $W$ and $V/W$ are in fact irreducible in the definable category, i.e., any proper definable submodule is trivial : let $W_1$ be a finite $R$-module contained in $W$ and let $0\neq w_1$; then $K_W\cdot w_1=W_1$ is infinite, contradiction (the proof for $V/W$ is similar). Till the end of the proof, we use irreducibility in this stronger sense.

 We may suppose that $\ann_R(W)=\ann_R(V/W)$. Indeed, let $0 \neq r \in \ann_R(W) - \ann(V/W)$. The definable subgroups $\im(r)$ and $W \subseteq \ker(r)$ are $R$-invariant since $R$ is commutative. Then $\ker(r) / W$ and $(\im(r) + W)/W$ are $R$-invariant, and by irreducibility of $V/W$, we have $V = \im(r) +W$ and $\ker(r)=W$. Moreover, $\dim(V)=\dim(\im(r))+\dim(\ker(r))=\dim(\im(r))+\dim(\ker(r))-\dim(W\cap \im(r))$; thus $\dim(W\cap \im(r))=0$ and $\im(r)\cap W$ is a proper finite $R$-module contained in W, hence trivial.
 \\
Conversely, let $0\neq r\in \ann_R(V/W)-\ann_R(W)$. Then $\im(r)=W$ by hypothesis and irreducibility of $W$. But $\ker(r)\cap W$ is a definable $R$-module properly contained in $W$ and so trivial by irreducibility of $W$. Therefore $V=W\oplus \ker(r)$.
\\
Now, let $x\in \ann_R(W)=\ann_R(V/W)=I$; since $R=R^p$, there exists $y\in R$ such that $x=y^p$. In the domain $R/I$, we have $\overline{y}^p = \overline{x} = 0$, so $y \in I$. Then $y \cdot V \subseteq W$ and $y^2 \cdot V = \{0\}$. Finally, $x = 0$ and $I=0$, so $R$ embeds in $K_W$ and in $K_{V/W}$. Therefore $R$ acts by automorphisms on $W$ and $V/W$. 
\\
Now, we may suppose that $R$ acts by automorphisms on $V$ itself. Indeed, let $0\neq r\in R$ and suppose $\ker(r)\neq 0$. Then by irreducibility of $W$, $\ker(r)\cap W=0$; so by irreducibility of $V/W$, $V=W\oplus \ker(r)$. 
\\
In $\End_{\text{def}}(V)$, we consider the subring $K=R\cdot R^{-1}$; this is an invariant and unbounded field. Let $Z\subseteq K$ be a type-definable subset. For $0\neq v\in V$, the evaluation map $K\rightarrow V, k\mapsto k\cdot v$ is injective since $K$ acts by automorphisms; thus $\dim(Z)=\dim(Z\cdot v)\leq \dim(V)$. By Lemma \ref{field definability}, $K$ is in fact a definable field. As a consequence, $W$ is a $K$-subspace of finite linear dimension and admits a definable $R$-invariant complement.
\end{proof}

\begin{corollary}\label{complement irreductible} (compare with \cite[Proposition 2.13]{Tind})
In a finite-dimensional$^{\circ}$ theory, let $T$ be a definable connected abelian group with trivial $p$-torsion and $V$ be a definable connected $T$-module that is $p$-elementary as an abelian group. Suppose that $V^{T} = 0$. Then every definable connected submodule of $V$ admits a definable connected $T$-invariant complement.
\end{corollary}

\begin{proof}
We call the following "Maschke property" :
\begin{quote}
Every definable connected submodule has a complement definable connected submodule. 
\end{quote}
First, we prove the Maschke property for $V$ of $T$-length two. Let $W$ be a proper definable connected submodule; then we may suppose that $W$ and $V/W$ are both $T$-irreducible.

Let $R$ be the subring of $\End_{\text{def}}(V)$ generated by $T$. Note that $R$ is commutative, invariant and unbounded since $T$ is definable, connected, and acts non-trivially on $V$. We have $R^p=R$ since $T$ is abelian and $T^p=T$. Moreover, $W$ and $V/W$ are $R$-irreducible; note also that $R/\ann_R(W)$ and $R/\ann_R(V/W)$ are unbounded since $T$ acts non-trivially on $W$ and $V/W$ by hypothesis and by Corollary \ref{facteur de composition}. Then we are done by Proposition \ref{anneau invariant}.  

By Corollary \ref{facteur de composition}, the hypothesis $V^{T}=0$ is still true for definable connected $T$-invariant subgroups and sections of $V$. We proceed by induction on the length of $V$ as a definable connected $T$-module to reduce to the case $\lg_T(V)=2$. Let $V$ be a definable connected $T$-module such that all definable connected $T$-modules of shorter length have Maschke property. Let $W$ be a definable connected $T$-submodule of $V$.
\begin{itemize}
\item
Suppose $\lg(W) > 1$. Then by definition, there is a definable connected $T$-module $0 < V_0 < W$; let $\pi_0$ stand for projection modulo $V_0$.

By induction, $\overline{V} = V/V_0$ has the Maschke property, so there is a definable connected $T$-module $\overline{X}$ with $\overline{V} = \overline{W} \oplus \overline{X}$. Let $X = (\pi_0^{-1}(\overline{X}))^o$. Then $V_0 \leq X < V$ is a definable connected $T$-submodule with $V = W + X$.

By induction, $X$ has the Maschke property, so there is a definable connected $T$-module $W'$ with $X = V_0 \oplus W'$. Hence on the one hand, $V = W + W'$, and on the other hand $\pi_0(W \cap W') \leq \overline{W} \cap \overline{X} = 0$, so $W \cap W' \leq V_0 \cap W' = 0$: thus $W'$ is a direct complement of $W$ in $V$.
\item
Suppose $\lg(V/W) > 1$. Then by definition there is a definable connected $T$-module $W < V_1 < V$; let $\pi_W$ stand for projection modulo $W$.

We may assume $W \neq 0$, so $\overline{V} = V/W$ has the Maschke property and there is a definable connected $T$-module $\overline{X}$ with $\overline{V} = \overline{V_1} \oplus \overline{X}$; let $X = (\pi_W^{-1}(\overline{X}))^o$; clearly $W \leq X < V$.

Since $V_1$ has the Maschke property, there is a definable connected $T$-module $Y \leq V_1$ with $V_1 = W \oplus Y$; since $X$ has too, there is $Z \leq X$ with $X = W \oplus Z$. Now let $W' = Y+Z$, again a definable connected $T$-module. Then on the one hand, $V = W + W'$, and on the other hand, if $y + z \in W$ with obvious notations, then $\pi_W(y) = - \pi_W(z) \in \overline{Y} \cap \overline{Z}\leq \overline{V_1}\oplus \overline{X} = 0$, so that $y \in Y \cap W = 0$ and $z = 0$ likewise; this shows $W \cap W' = 0$. Thus $W'$ is a direct complement of $W$ in $V$.
\end{itemize}

\end{proof}
\begin{corollary}
 In a finite-dimensional$^{\circ}$ theory, let $T$ be a definable connected abelian group with trivial $p$-torsion and $V$ be a definable connected $T$-module that is $p$-elementary as an abelian group. Suppose that $V^{T} = 0$. Then $V=V_1\oplus\dots\oplus V_r$, where each $V_i$ is a definable connected irreducible $T$-module.
\end{corollary}
Note that we do not consider the nilpotent case since  we need to work with a \textit{commutative} ring $R$ such that $R^p=R$ in order to implement our strategy.
\section{Lie rings}
\subsection{Cohomological results}
In \cite{HS}, Hochschild and Serre define a cohomology theory for Lie algebras that parallels the theory for groups. The definition can easily be adapted to the more general context of Lie rings; nevertheless, general cohomological results for Lie algebras may fail when passing to Lie rings. Recall that a Lie ring $(\mathfrak{g},+,[,])$ is an abelian group equipped with a bi-additive and anti-symmetric function, $[,]$, that satisfies the Jacobi identity.
\begin{definition}
A $\mathfrak{g}$-module is the data of ($\mathfrak{g}$,$A$), where $\mathfrak{g}$ is a Lie ring and $A$ is an abelian group, equipped with a Lie ring action "$\cdot$" of $\mathfrak{g}$ on $A$. Equivalently, this amounts to giving a Lie ring homomorphism $\rho$ from $\mathfrak{g}$ to $\End(A)$.
\end{definition}
Recall that a $\mathbb{Z}$-multilinear map $f$ is alternating if $f(x_1,\dots,x_n)=0$ whenever two arguments among the $x_i$ are equal.
\begin{definition}
Let $(\mathfrak{g},A)$ be a $\mathfrak{g}$-module. We define a cochain complex as follows:
\begin{itemize}
    \item $C^0(\mathfrak{g},A)=A$
    \item For every integer $n\geq 1$, $C^{n}(\mathfrak{g},A)=\{f\in \Hom(\mathfrak{g}^n,A) : f \text{ is alternating}\}$, where $\Hom(\mathfrak{g}^n,A)$ denotes the set of maps from $\mathfrak{g}^n$ to $A$, which are $\mathbb{Z}$-multilinear.  
   
\end{itemize}
The differential $d_n$ from $C^n(\mathfrak{g},A)$ to $C^{n+1}(\mathfrak{g},A)$ is defined, for every map $f$ in $C^n(\mathfrak{g},A)$, by:

\[\begin{aligned}& d_nf(g_1,\dots,g_{n+1})\\
& =\sum_{1\leq s < t\leq n+1}(-1)^{s+t-1}f([g_s,g_t],g_1,\dots, \hat{g_s},\dots,\hat{g_t},\dots, g_{n+1})\\
& + \sum_{1\leq s \leq n+1}(-1)^sg_s\cdot f(g_1,\dots,\hat{g_s},\dots,g_{n+1})\end{aligned}\]
 \end{definition}
  The Lie ring $\mathfrak{g}$ acts on the cochain groups in a way that is compatible with the differentials: 
 
 For all $f\in C^n(\mathfrak{h},A)$, where $\mathfrak{h}$ is an ideal, and for all $x\in \mathfrak{g}$; we define $x\cdot f$ by the formula \[(x\cdot f)(y_1,\dots,y_n)=x\cdot f(y_1,\dots,y_n)-\sum_{1\leq i \leq n}f(y_1,\dots,y_{i-1},[x,y_i],y_{i+1},\dots,y_n),\] for all $y_1,\dots, y_n\in \mathfrak{\mathfrak{h}}$.
 \\
 \\
 For $n=0$ and $n=1$, it is also possible to give a convenient algebraic description of the cohomology groups.
 \begin{fact}
 \begin{enumerate}
     \item $H^0(\mathfrak{g},A)=A^{\mathfrak{g}}=\{a\in A : g\cdot a=0, \text{for all $g\in \mathfrak{g}$}\}$.
     \item $H^1(\mathfrak{g},A)=\Der(\mathfrak{g},A)/\IDer(\mathfrak{g},A)$, where \[\Der(\mathfrak{g},A)=\{f\in \Hom(\mathfrak{g},A)  : f([x,y])=x.f(y)-y.f(x)\}\] and \[\IDer(\mathfrak{g},A)=\{f\in \Hom(\mathfrak{g},A) : f(x)=x.a~ \text{for some $a$ in $A$}\}.\] 
 \end{enumerate}
     
 \end{fact}
All the statements that we proved regarding the cohomology of groups are still true when passing to Lie rings, and the proofs require minor adaptations. In particular, we have the following:
\begin{fact}\label{inflation-restriction lie algebras}
    \begin{itemize}
        \item Every short exact sequence of $\mathfrak{g}$-modules $0\rightarrow A\overset{i}{\rightarrow}B\overset{p}{\rightarrow}C\rightarrow0$ induces a "long" exact sequence of cohomology groups of the following form : \[\begin{aligned}&A^{\mathfrak{g}}\rightarrow B^{\mathfrak{g}}\rightarrow C^{\mathfrak{g}}\\
    &\overset{\delta_0}{\rightarrow}H^{1}(\mathfrak{g},A)\rightarrow H^{1}(\mathfrak{g},B)\rightarrow\ H^{1}(\mathfrak{g},C)\end{aligned}\]
        \item For an ideal $\mathfrak{h}$, the sequence \[0\rightarrow H^1(\mathfrak{g}/\mathfrak{h},A^{\mathfrak{h}})\overset{\inf}{\rightarrow} H^1(\mathfrak{g},A)\overset{\res}{\rightarrow} H^1(\mathfrak{h},A)^{\mathfrak{g}/\mathfrak{h}}\] is exact.  
    \end{itemize}
\end{fact}

\begin{Lemma}\label{fidélité}
    Let $(\mathfrak{g},A)$ be a $\mathfrak{g}$-module and let $\mathfrak{h}$ be a central ideal of $\mathfrak{g}$, such that $\mathfrak{h}\leq C_{\mathfrak{g}}(A)$. Suppose that $A^{\frak{g}}=0$. Then $H^1(\mathfrak{g}/\mathfrak{h}, A)\simeq H^1(\mathfrak{g}, A)$.
\end{Lemma}
\begin{proof}
Since $A^{\frak{h}}=A$, by Fact \ref{inflation-restriction lie algebras}, the sequence $0\rightarrow H^1(\mathfrak{g}/\mathfrak{h},A)\overset{\inf}{\rightarrow} H^1(\mathfrak{g},A)\overset{\res}{\rightarrow} H^1(\mathfrak{h},A)^{\frak{g}/\frak{h}}$ is exact. Therefore, it suffices to prove that $H^1(\mathfrak{h},A)^{\mathfrak{g}/\mathfrak{h}}=0$.

Let $f\in \Der(\mathfrak{h},A)$; suppose that for all $x\in \mathfrak{g}$, $[x\cdot f]=[0]$. For all $x\in \mathfrak{g}$ there exists an $a_x\in A$ such that for all $h\in \mathfrak{h}$ one has $x\cdot f(h)-f([x,h])=x\cdot f(h)=h\cdot a_x=0$ since $\mathfrak{h}$ acts trivially and $\mathfrak{h}\leq Z(\mathfrak{g})$. In other words, for all $x\in \mathfrak{g}$ and $h\in \mathfrak{h}$, the equality $x\cdot f(h)=0$ holds, and thus $f(h)\in A^\mathfrak{g}=0$; finally, we have $f=0$. 
\end{proof}
There is also a minor variation concerning the abelian case for irreducible module (Lemma \ref{abelian case finite-dimensional}).
\begin{Lemma}
In a finite-dimensional$^{\circ}$ theory, let $\mathfrak{g}$ be a definable connected abelian Lie ring and $A$ be a definable connected $\mathfrak{g}$-module. Suppose that $A$ is $\mathfrak{g}$-irreducible and $A^\mathfrak{g}=0$. Then $H^1(\mathfrak{g},A)=0$.
\end{Lemma}
\begin{proof}
    First, we way suppose that $\mathfrak{g}$ acts faithfully by Lemma \ref{fidélité} and so we have an embedding $\rho :  \mathfrak{g} \hookrightarrow \End_{\text{def}}(A) $. Let $R$ be the subring of $\End_{\text{def}}(A)$ generated by $\mathfrak{g}$; by Fact \ref{linearisation} , $C_{\End_{\text{def}}(A)}(R)=\mathbb{K}$ is a definable skew-field containing $R$, in particular $\mathfrak{g} \le R$ acts by automorphisms. Up to considering its center, we may suppose that $\mathbb{K}$ is commutative.
    \\
    Now, let $\delta : \mathfrak{g}\rightarrow A$ be a $\mathfrak{g}$-derivation. Let $y\in \mathfrak{g}-\{0\}$; then $\rho(y)$ is invertible. Let $a=\rho(y)^{-1}\cdot\delta(y)$. For $x\in \mathfrak{g}$, by commutativity and the definition of a $\mathfrak{g}$-derivation, we have :
    \[\rho(x)\cdot a=\rho(x)\cdot \rho(y)^{-1}\cdot \delta(y)=\rho(y)^{-1}\cdot \rho(x)\cdot \delta(y)=\rho(y)^{-1}\cdot \rho(y)\cdot \delta(x)=\delta(x). \]
    Thus $\delta(x)=\delta_a(x)$, where $\delta_a(x)=\rho(x)\cdot a$ is the $\mathfrak{g}$-inner derivation associated with $a$. This shows that $H^1(\mathfrak{g},A)=0$.
\end{proof}
Finally, similar to Theorem \ref{g-module} and Corollary \ref{facteur de composition}, we have the following:
\begin{theorem}\label{cohomology Lie rings}
  \begin{itemize}
    \item In a finite-dimensional$^{\circ}$ theory, let $\mathfrak{g}$ be a definable connected nilpotent Lie ring and $A$ be a definable connected $\mathfrak{g}$-module. Suppose that $A^{\mathfrak{g}}=0$. Then $H^1(\mathfrak{g},A)=0$.
    \item In a finite-dimensional$^{\circ}$ theory, let $\mathfrak{g}$ be a definable connected nilpotent Lie ring and $A$ be a definable connected $\mathfrak{g}$-module. Suppose that $A^{\mathfrak{g}}=0$. Then for all definable connected submodules $V\leq U$ contained in $A$, one has $(U/V)^{\mathfrak{g}}=0$.
\end{itemize}  
\end{theorem}

\subsection{Definable Maschke's theorem for restricted Lie algebras}
Regarding Maschke's theorem, we must pass to Lie algebras since the abstract Lie ring structure seems to be insufficient to capture a good notion of $p$-divisibility/$p$-unipotency. First, we clarify what semi-simplicity means for Lie algebras in positive characteristic. We start with the algebra of endomorphisms of a finite-dimensional vector space, for which linear algebra provides a natural notion of semi-simple element. For the following developments, we refer to Chapter II of \cite{Stra}.
\\
Let $V$ be a finite-dimensional vector space over a field $\mathbb{K}$ of characteristic $p>0$; consider the associative $\mathbb{K}$-algebra $\mathfrak{g} = \End(V)$. We aim to characterize semi-simple endomorphisms purely in terms of algebraic operations in $\mathfrak{g}$ without reference to the minimal polynomial. To this end, raising to the $p$-th power will play a key role; we begin with the following definition:

\begin{definition}
An element $x$ of $\mathfrak{g}$ is \emph{algebraically semi-simple} if 
\[
x \in \langle \bigcup_{i\geq 1} x^{p^i} \rangle_{vect}.
\]
\end{definition}

The following fact shows that the two notions of semi-simplicity coincide:

\begin{fact}\cite[Proposition 3.3, chap. 2]{Stra}
An element $x$ of $\mathfrak{g}$ is algebraically semi-simple if and only if it induces a semi-simple endomorphism.
\end{fact}
For a Lie algebra over a field of characteristic $p>0$, one needs an analogue of the $p$-power map; this leads to the notion of \emph{restricted Lie algebras}. They naturally appear in the classification of simple Lie algebras of finite linear dimension over an algebraically closed field of positive characteristic. For matrix Lie algebras, the $p$-function is really the $p$-power map. Note that, in general, owing to the non-commutativity encoded in the Lie bracket, the $p$-map is not additive.
\begin{definition}
A Lie algebra $\mathfrak{g}$ over a field $\mathbb{K}$ of characteristic $p>0$ is \emph{restricted} if there exists a map $[p]:\mathfrak{g} \to \mathfrak{g}$, $x \mapsto x^{[p]}$, satisfying the following conditions:
\begin{enumerate}
    \item $[x^{[p]},y] = \ad_x^p(y)$ for all $y \in \mathfrak{g}$.
    \item $(\lambda x)^{[p]} = \lambda^p x^{[p]}$ for all $\lambda \in \mathbb{K}$.
    \item $(x+y)^{[p]} = x^{[p]} + y^{[p]} + \sum_{1\le i \le p-1} s_i(x,y)$, where $\ad_{x\otimes X + y\otimes 1}^{p-1}(x \otimes 1) = \sum_{1\le i \le p-1} i s_i(x,y) \otimes X^{i-1}$ in $\mathfrak{g} \otimes_\mathbb{K} \mathbb{K}[X]$, for all $x,y \in \mathfrak{g}$.
\end{enumerate}
Moreover, if $\mathfrak{g}$ is abelian, then $s_i(x,y) = 0$ for all $i \in \{1,\dots,p-1\}$ and all $x,y \in \mathfrak{g}$.
\end{definition}

Note that an associative $\mathbb{K}$-algebra with its natural Lie algebra structure is restricted when $[p]$ is taken as the $p$-th power. In particular, the algebra of endomorphisms of a vector space over a field of positive characteristic and the universal enveloping algebra of a Lie algebra over a field of positive characteristic have a natural structure of restricted Lie algebras.

\begin{definition}
Let $(\mathfrak{g},[p])$ be a restricted Lie algebra over a field of characteristic $p>0$. An element $x$ is \emph{semi-simple} if 
\[
x \in \langle \bigcup_{i \ge 1} x^{[p]^i} \rangle_{vect}.
\]

In a finite-dimensional theory, a definable connected restricted Lie algebra $(\mathfrak{t},[p])$ is a \emph{torus} if it is abelian and $ker([p]) = 0$.
\end{definition}

\begin{Remark}
In the linear finite-dimensional case, an element of a restricted Lie algebra over a perfect field of characteristic $p>0$ is semi-simple if and only if it belongs to an abelian subalgebra stable under $[p]$, such that $ker([p])=0$ in this subalgebra (see \cite[Proposition 3.3 and Theorem 3.7]{Stra}). 
\end{Remark}
Now, we may use Proposition \ref{anneau invariant} and argue as in the group case to prove the following:
\begin{corollary}
In a finite-dimensional$^{\circ}$ theory, let $\mathfrak{t}$ be a torus and $V$ a definable connected $\mathfrak{t}$-module. Assume that $V$ is abelian $p$-elementary and $\rho(t^{[p]})(x) = \rho(t)^p(x)$ for all $t \in \mathfrak{t}$ and $x \in V$, where $\rho: \mathfrak{t} \to \End(V)$ is the map given by the action. Suppose that $V^{\mathfrak{t}} = 0$. Then $V=V_1\oplus\dots\oplus V_r$, where each $V_i$ is a definable connected irreducible $\mathfrak{t}$-module.
\end{corollary}
\subsection{Frattini's argument for Lie rings}
In the theory of finite groups, the $p$-Sylow subgroups play a fundamental role and often allow to obtain structural results about the ambient group. In particular, Frattini's argument states that, for a finite group $G$, a normal subgroup $H$, and a $p$-Sylow subgroup $S$ of $H$, we have the following :  
\[
G = H  N_G(S).
\]  
The argument is based on conjugacy results for $p$-Sylow subgroups and thus generalises poorly to the setting of Lie rings: even when one has a suitable analogue of $p$-Sylow subgroups, the "inner" maps are derivations, therefore not automorphisms, and passing through the exponential map, when it is well-defined, only partially remedies this defect. However, the conclusion of Frattini's argument remains valid for finite-dimensional (in the usual linear sense) Lie algebras when $p$-Sylow subgroups are replaced by Cartan subalgebras. The proof relies on the cohomological study of nilpotent Lie algebras (see \cite[Theorem 2.1]{Barn2}).  

In light of our cohomological results and following the argument of \cite[Theorem 2.1]{Barn2}, we now prove an analogous result in the context of finite-dimensional$^{\circ}$ theories.

\begin{definition}
In a finite-dimensional$^{\circ}$ theory, let $\mathfrak{g}$ be a definable connected Lie ring and let $\mathfrak{c}$ be a definable connected nilpotent subring such that  
\[
N_{\mathfrak{g}}^{\circ}(\mathfrak{c}) = \mathfrak{c}.
\]  
We then say that $\mathfrak{c}$ is a \emph{Cartan subring}.
\end{definition}

\begin{proposition}
In a finite-dimensional$^{\circ}$ theory, let $\mathfrak{g}$ be a definable connected Lie ring, and let $\mathfrak{c}$ be a definable connected Cartan subring of a definable connected ideal $\mathfrak{i}$.  Then  
\[
\mathfrak{g} = \mathfrak{i} + N^{\circ}_{\mathfrak{g}}(\mathfrak{c}).
\]
\end{proposition}

\begin{proof}
We set $\mathfrak{n} = N_{\mathfrak{g}}(\mathfrak{c})$ and we consider the definable connected $\mathfrak{c}$-module $A = \mathfrak{g} / \mathfrak{n}$.

Let us prove that $A^{\mathfrak{c}} = 0$ for the adjoint action.  Indeed, let $\overline{a}\in A^\mathfrak{c}$. Lift $\overline{a}\in A$ to $a \in \mathfrak{g}$. Then $[\mathfrak{c},a]$ is a connected subgroup of $\mathfrak{i} \cap \mathfrak{n}$, hence contained in $N_{\mathfrak{i}}^{\circ}(\mathfrak{c})=\mathfrak{c}$. Thus $a \in \mathfrak{n}$ and $\overline{a}=0$.

Now, consider the $\mathfrak{c}$-module
\[
B=(\mathfrak{g}/\mathfrak{n}) / ((\mathfrak{i} + \mathfrak{n}) / \mathfrak{n}) \simeq \mathfrak{g} / (\mathfrak{i} + \mathfrak{n}).
\]  
Since $B$ is a quotient of $A$,  $B^\mathfrak{c}=0$ by Theorem \ref{cohomology Lie rings}. But $\mathfrak{c} \le \mathfrak{i} \triangleleft \mathfrak{g}$, so $[\mathfrak{g},\mathfrak{c}]\le \mathfrak{i} \le \mathfrak{i}+\mathfrak{n}$, proving $B = B^\mathfrak{c} = 0$. Therefore $\mathfrak{g}=\mathfrak{i}+\mathfrak{n}$. By connectedness, $\mathfrak{g} = \mathfrak{i}+\mathfrak{n}^{\circ}$.
\end{proof}
Some applications of Frattini's argument will be given in \cite{TZ} (in collaboration with Tindzhogo Ntsiri). Notably, we prove that the Frattini subring is nilpotent in the context of connected soluble Lie rings of finite Morley rank.

A generalisation of some parts of this paper is currently under investigation by Moreno Invitti.

\bibliography{cohomology}
\bibliographystyle{plain}

Samuel Zamour

Université Paris-Saclay

Laboratoire de Mathématiques d'Orsay (LMO), équipe "Arithmétique et Géométrie algébrique" 

e-mail: samuel.zamour@u-pec.fr
\end{document}